\documentclass[11pt]{article}
\usepackage{amssymb,amsmath}

\newtheorem{theorem}{Theorem}[section]
\newtheorem{lemma}[theorem]{Lemma}

\newtheorem{proposition}[theorem]{Proposition}

\newcommand{\g}{\mathfrak{g}}
\newcommand{\h}{\mathfrak{h}}

\begin{document}

\title{Chern classes  of quantizable coisotropic bundles}

\author{Vladimir Baranovsky}
\date{April 1, 2021}

\maketitle

\abstract{Let $M$ be a smooth algebraic variety of dimension $2(p+q)$
 with an algebraic symplectic form 
and a compatible deformation quantization $\mathcal{O}_h$ of the structure sheaf.
Consider a smooth coisotropic
subvariety $j: Y \to M$ of codimension $q$ and a vector bundle $E$ on $Y$.
We show that if $j_* E$ admits a deformation quantization (as a module) 
then its characteristic class $\widehat{A}(M) exp(-c(\mathcal{O}_h)) ch(j_* E)$ lifts
 to a cohomology group associated to the null foliation 
of $Y$. Moreover, it can only be
nonzero in degrees $2q, \ldots, 2(p+q)$.  
For Lagrangian $Y$ this reduces to a single degree $2q$. Similar results hold
in the holomorphic category.

\section{Introduction}

Let $(M, \omega)$ be a smooth variety of dimension $2(p+q)$
over a field $k$ of 
characteristic zero, with an algebraic symplectic form $\omega$ 
(or corresponding holomorphic objects over $k = \mathbb{C}$). 
We assume that the structure sheaf $\mathcal{O}_M$ admits
a compatible deformation quantization  $\mathcal{O}_h$ and fix a choice
of such quantization.   In other words, $\mathcal{O}_h$ is
a sheaf (in the Zariski or analytic topology, respectively) of complete, 
separated and flat $k[[h]]$ algebras such that 
$\mathcal{O}_h/ h \mathcal{O}_h \simeq \mathcal{O}_M$ and 
if $a \mapsto a_0$ is the quotient map then 
$$
a * b - b * a = h P(d a_0, d b_0) (mod\  h^2)
$$
where $P \in H^0(M, \Lambda^2 T_M)$ is the Poisson bivector corresponding
to the symplectic form $\omega$ under the isomorphism $T_M \to \Omega^1_M$, 
$v \mapsto \iota_v (\omega)$. 

\bigskip
\noindent
Consider a coherent sheaf $E_h$ of $\mathcal{O}_h$-modules which is
complete, separated and flat over $k[[h]]$. See \cite{KS} for a general overview
of modules over deformation quantization. 
We view $E_h$ as a quantization of its
``principal symbol''  $\sigma(E_h) = E_h/h E_h$, a coherent sheaf of 
$\mathcal{O}_M$-modules. A broad, but difficult, question is to establish necessary 
and sufficient
conditions  which would imply existence of $E_h$. One way to 
simplify the situation is to assume that the support $Y \stackrel{j} \rightarrow M$ of $E_h$ is smooth, and
in fact $E_h/h E_h$ is a direct image $j_* E$ 
of a locally free sheaf $E$ on $Y$.

\bigskip
\noindent
A straigtforward observation, which we recall below, is that in this case $Y$ should be 
coisotropic, i.e. if 
 $N$ is the normal bundle of $Y$ in $M$ then the projection of $P$
to $H^0(Y, \Lambda^2 N)$ should be zero. 
Then $p = \frac{1}{2}(dim\;M - 2\;codim\;Y)$ is a non-negative integer. 
When $p=0$, i.e. $Y$ is Lagrangian, papers 
\cite{BGKP} and \cite{BC} establish
necessary and sufficient conditions for existence of $E_h$. First, the associated
projective bundle $\mathbb{P}(E)$ on $Y$ should admit a flat 
algebraic connection. In particular,  the Chern character of $E$ equals 
$e \cdot exp(c_1(E))$ with $e = rk\;E$. To formulate the remaining conditions, 
recall that a choice of $\mathcal{O}_h$ induces the  
Deligne-Fedosov class $c(\mathcal{O}_h) \in 
\frac{1}{h} H^2_{DR}(X)[[h]]$ of the form 
$$
c(\mathcal{O}_h) = \frac{1}{h} [\omega] + \omega_0 + 
h \omega_1 + h^2 \omega_2 + \ldots
$$
Note that our present indexing of the coefficients is shifted by 1 as compared
to that of \cite{BC}.

The Lagrangian property of $Y$ implies that $[\omega]$ restricts to 
zero on $Y$. In \cite{BGKP}, \cite{BC} it was shown that existence of
quantization also implies $\omega_i|_Y = 0$, for the cohomology classes
$\omega_i$ with  $i \geq 1$. This 
may be viewed as a strengthened Lagrangian condition which depends
on the choice of $\mathcal{O}_h$. 
As for the cohomology class 
$\omega_0$, it admits a canonical lift
(depending on $\mathcal{O}_h$)  to $H^2_{DR}(Y, \Omega^{\geq 1})$ which
is involved in the equation:
\begin{equation}
	\label{lie-class}
\frac{1}{e} c_1(E) = \omega_0 + \frac{1}{2} c_1(K_Y)
\end{equation} 
where $K_Y$ is the canonical bundle of $Y$. 
Square root of the canonical class has appeared in the contact setting in 
$\cite{Ka}$. In the holomorphic setting, quantization of the square root of
the canonical class is due to D'Agnolo and Schapira, see \cite{DS}. 
See also \cite{Bo} and \cite{NT2}. 

\bigskip
\noindent
Returning to the case of a non-necessarily Lagrangian smooth coisotropic $Y$, 
define the  ``quantum Chern character" class 
 considered in \cite{NT} (see also \cite{BNT} and \cite{CFW} for  later 
 results):
$$
\tau(E_h) : = \widehat{A}(T_M) exp(-c(\mathcal{O}_h)) ch(\sigma(E_h)).
$$
where the $\widehat{A}$-genus is recalled in Section 2.2.
We understand $\tau(E_h)$ as a class with values in the de Rham cohomology 
$H^*_{DR, Y} (M)((h))$ with support at $Y$, which can be
identified with the de Rham cohomology of $Y$, due to smoothness. 
The purpose of this note is the following result
\begin{theorem} 
	\label{main}
	If the principal symbol sheaf $\sigma(E_h)$ is isomorphic to 
the 
direct image $j_* E$ of a locally free sheaf $E$ on a smooth coisotropic
subvariety $j: Y \to M$ of codimension $q$ then the class 
$\tau(E_h) \in H^*_{DR, Y}(M)((h))$ is zero except 
in degrees $2q, \ldots, 2(p+q) = \dim_k M $. 

Moreover, if $\Omega_{\mathcal{F}} \subset \Omega^1_Y$ is the sheaf of 1-forms
that vanish on the null-foliation (or characteristic foliation)  $\mathcal{F} \subset T_Y$ and
$F^r \Omega^\bullet_Y$ is the ideal in the de Rham complex generated by 
the $r$-th power of $\Omega_{\mathcal{F}}$, then $\tau(E_h)$ is in the image of 
the map 
$$
\bigoplus_{r \geq 0}^p H^{2r}(Y, F^r \Omega^\bullet_Y)  ((h))
\to   \bigoplus_{r \geq 0}^p H^{2r}(Y, \Omega^\bullet_Y)  ((h))
\simeq  \bigoplus_{r \geq 0}^p H^{2r+2q}_Y(M,  \Omega^\bullet_M)  ((h))
$$
\end{theorem}
Our strategy is an application of formal geometry 
and the Gelfand-Fuks map: first use Riemann-Roch 
theorem to replace $\tau(E_h)$ by an element $\tau_Y(E) \in H^\bullet_{DR}(Y)((h))$; 
then show that after completion both the quantized functions and the quantized
module are isomorphic to standard objects and construct a Harish-Chandra
torsor (foliated over $\mathcal{F}$) 
and a Lie algebra cohomology class that induces $\tau (E_h)$ via
Gelfand-Fuks map. At this point the vanishing reduces to a 
vanishing in Lie algebra cohomology for which we use a Lie Algebraic version of
the Index
Theorem, cf. \cite{NT}, \cite{BNT}, \cite{CFW}, \cite{DC}, \cite{PPT}, \cite{GKN},
\cite{GLX}, 
 and that fact that in the symplectic situation the trace map can be
defined on negative cyclic homology. 

Alternatively, one could use an observation due to B. Tsygan that the
Chern character of a perfect complex factors through negative cyclic homology of
its derived endomorphism algebra, but in the algebraic geometry setting
the Lie algebra cohomology route seems a bit shorter. 

\bigskip
\noindent 
\textbf{Remarks.} (i) In a forthcoming paper with V. Ginzburg, cf. 
\cite{BG},  we prove a similar
statement for quantizable sheaves with arbitrary support, including the fact that $\tau(E_h)$ 
agrees with the general Connes' Chern character and that the 
Algebraic Index Theorem holds for general algebraic varieties. Since in general a 
formal completion of a quantized sheaf will not be isomorphic to a ``standard 
formal model", methods of formal geometry do not apply for general sheaves. 

(ii) It would be very interesting to relate our main theorem to Bordemann's criterion 
for existence of second order quantization (modulo $h^3$). This might depend on 
what can be said about the Atiyah-Molino class of the charactersitic foliation of
$Y$. See \cite{Bo} for more details.

\bigskip
\noindent 
The paper is organized as follows. 
In Section 2 we recall standard constructions related to foliations, characteristic classes
and use Riemann-Roch theorem to reduce the main result to a cohomology class on $Y$.
In Section 3 we recall definitions related to Harish Chandra pairs and torsors and the
Gelfand Fuks map. We further state the Lie cohomology Algebraic Index Theorem and prove a vanishing
result for the class involved. The conceptual reason for the vanishing is that Connes Chern
character with values in periodic cyclic homology lifts to negative cyclic homology.  
In Section 4 we prove the main result by constructing two Harish-Chandra torsors that
induce the class under consideration, and then invoking the vanishing of Section 3. 

\bigskip
\noindent
\textbf{Acknowledgements.} The author thanks V. Ginzburg,
 J. Pecharich and T. Chen for 
the useful conversations. This work was supported by 
the Simons Collaboration Grant \#281515. 

\section{Preliminaries and notation.}

\subsection{Null foliation and a filtration on the de Rham complex.}

We start by assuming that a pair $(\mathcal{O}_h, E_h)$ is given 
as in the introduction and that $\sigma(E_h) = E_h/h E_h$ is the 
direct image $j_* E$ of a locally free sheaf supported on a smooth 
subvariety $Y$. We use the same notation $E$ for associated vector bundle
on $Y$. 

If $I \subset \mathcal{O}_M$ is the ideal sheaf of functions vanishing on 
$Y$ and $x \in \mathcal{O}_h$ is a local section projecting to $I$, then 
$x \cdot E_h \subset h E_h$. If $y$ is another such section, it follows that
the commutator $(x \cdot y - y \cdot x)$ sends $E_h$ to $h^2 E_h$ and 
hence the image of $\frac{1}{h}(x \cdot y - y \cdot x)$ in $\mathcal{O}_M$
also annihilates $j_*E$, i.e belongs to $I$. 
Thus, the ideal sheaf $I$ is closed with respect 
to the Poisson bracket induced by the Poisson bivector $P$. 

If $N$ is the normal bundle of $Y$ in $M$, we can restate this by saying that
$P|_Y$ projects to the zero section in $H^0(Y, \Lambda^2 N)$, and then 
the same restriction defines a section in $H^0(Y, T_Y \otimes N)$.
We can view the latter as a morphism $N^\vee \to T_Y$ and it is easy
to check that it is an embedding of vector bundles. Using $j_* N^{\vee}
\simeq I/I^2$ we can write an explicit local formula for it:
$$
N^\vee \ni x \mapsto P|_Y(dx, \cdot) \in T_Y.
$$
Denote by $\mathcal{F} \subset T_Y$ the image of this embedding, i.e. 
the \textit{null-foliation} of $Y$. By the above, this sub-bundle
is \textit{involutive}, i.e.  closed with respect to the bracket of vector fields on $Y$
(since the Poisson bracket on $I/I^2$ is compatible with the bracket on vector fields).

The involutive property can be restated as follows. Let
$\Omega_{\mathcal{F}} = (T_Y/\mathcal{F})^\vee \subset \Omega^1_Y$
be the sheaf of 1-forms vanishing along $\mathcal{F}$ and denote by 
$F^1 \Omega^\bullet \subset \Omega^\bullet_Y$ the graded ideal
generated by $\Omega_{\mathcal{F}}$ in the sheaf of differential forms on $Y$,
viewed as a sheaf of graded commutative algebras. By a straightforward application
of the formula
\begin{equation}
	\label{differential}
d \omega(v_0, v_1) = v_0 \omega(v_1) - v_1 \omega(v_0) - \omega([v_0, v_1])
\end{equation}
 the involutive property of 
$\mathcal{F}$ is equivalent to the statement that $F^1 \Omega^\bullet$ is
a subcomplex of the de Rham complex. It follows that each power
of the ideal 
$F^k \Omega^\bullet := (F^1 \Omega^\bullet)^k$ is also a subcomplex.

The main message in this paper is that characteristic classes of interest
lift to the cohomology groups $H^{2r}(Y, F^r \Omega^\bullet)$. Note 
that the rank of $\Omega_{\mathcal{F}}$ is $2p = dim\;Y - rk(\mathcal{F}) =
dim\;M -2  \; codim\;Y$ hence for $r > 2p$ the relevant cohomology group 
vanishes as $F^r \Omega^\bullet$ is the zero subcomplex.  Moreover, the
particular class $\tau(E_h)$ vanishes for $r > p$. 

\subsection{Riemann-Roch Theorem and reduction to a class on $Y$.}
Recall that for a power series $G(z)$ with constant term $1$ and
a vector bundle $V$ of rank $v$ we can define its multiplicative
$G$-genus as a product $\prod G(z_i)$ where $z_1, \ldots z_v$ are
the Chern roots of $V$, i.e. formal variables such that the $l$-th
elementary symmetric function in $z_i$ is equal to the $l$-th Chern
class $c_l(V)$. Recall also that the $\widehat{A}$-genus 
$\widehat{A}(V)$ and the Todd genus $Td(V)$ correspond to  
$$
G_1(z) = \sqrt{\frac{z/2}{sinh(z/2)}} = \sqrt{\frac{z \cdot exp(-z/2)}{1 - exp(-z)}} 
\textrm{\ \   and\ \   } G_2(z) = \frac{z}{1- exp(-z)}
$$
respectively. 
The Todd genus is involved in the Grothendieck-Riemann-Roch theorem for
a closed embedding $j: Y \to M$ and a coherent sheaf $E$ on $Y$ 
(see Section 15.2 of \cite{Fu}): 
$$
ch(j_* E) = j_*[ch(E) Td(N)^{-1}],
$$
where $N$ is the normal bundle to $Y$. 
We use this formula to study the ``quantum"  class 
$$
\tau(E_h) = \widehat{A}(T_M) exp(-c(\mathcal{O}_h) ) ch(j_* E)
$$
which appears in the Index Theorem of \cite{BNT}, \cite{NT} and the Local Index Formula of \cite{DC}.
By Riemann-Roch and the Projection Formula we can rewrite this expression as
$$
\tau(E_h) = j_* \Big[\frac{\widehat{A}(T_M|_Y)}{Td(N)} 
exp(-c(\mathcal{O}_h|_Y) )ch(E) \Big]
$$
Note that $G_1(z)$ is an even function of $z$ hence $\widehat{A}(N)
= \widehat{A}(N^\vee)$. Using the multiplicative property of genus and 
the short exact sequences
$$
0 \to T_Y \to T_M|_Y \to N \to 0 \qquad 0 \to N^\vee \simeq \mathcal{F} \to T_Y \to Q \to 0
$$
(where the second short exact sequence is the definition of $Q$) we conclude that 
$$
\frac{\widehat{A}(T_M|_Y)}{Td(N)}  = \frac{\widehat{A}(Q) \cdot (\widehat{A}(N))^2}{Td(N)} = \widehat{A}(Q) exp(- \frac{c_1(N)}{2}) 
$$
Hence Theorem 1.1 reduces to the statement that the class 
$$
\tau_Y(E) = \widehat{A}(Q) exp(- \frac{c_1(N)}{2})  exp(-c(\mathcal{O}_h)|_Y) ch(E) 
$$
is in the image of $\bigoplus H^{2r}(Y, F^r \Omega^\bullet) \to \bigoplus  H^{2r}(Y, \Omega^\bullet) $ and vanishes in for $r > p$. 
Note that  when $M = Y$ and $M$ is projective the cohomology groups will be nonzero up
to degree $4p = 2 \dim M$.

\section{Lie Algebra Cohomology and Gelfand-Fuks map.}

The purpose of this section is to review some known results involving Lie
Algebra Cohomology and characteristic classes, as they apply to deformation
quantization; and also to fix the notation. 

\subsection{Lie algebra cohomology and a version of the Chern-Weyl map.}

We follow Section 10.1 in \cite{Lo} and Chapter 1.3 in \cite{Fu1}.
For a Lie algebra $\g$ over a field $k$ and a $\g$-module $V$ the Lie algebra
cohomology groups are defined as 
$
H^\bullet (\g; V) = Ext^\bullet_{U(\g)}(k, V)
$
where $U(\g)$ is the universal enveloping algebra. Using the standard 
Koszul resolution $\Lambda^\bullet(\g) \otimes U(\g) \to k$ over $U(\g)$, these can be
computed using  
 Chevalley-Eilenberg cochain complex $C^\bullet(\g; V)$:
 $$
 \ldots \to Hom_k(\Lambda^n \g, V) \to Hom_k(\Lambda^{n+1}(\g), V)  \to \ldots
 $$
 with a differential $d_{Lie}$, cf.
 10.1.6 in \cite{Lo}:  
 $$
 d_{Lie} \alpha(g_0, \ldots,g_n) = \hspace{4in} 
 $$
 \begin{equation}
 	\label{differential2}
 	= \sum_{i=0}^n (-1)^i g_i  \alpha(g_0, \ldots, \widehat{g}_i, \ldots, g_n) 
 + \sum_{i < j} (-1)^{i+j}\alpha([g_i, g_j], g_0, \ldots, \widehat{g}_i, \ldots,  \widehat{g}_j, \ldots, g_n)  
\end{equation}
  We view the elements 
 $\alpha \in Hom_k (\Lambda^\bullet(g), V)$ as skew-symmetric functions on $\g$ 
 with values in $V$. 
 
 \bigskip
 \noindent 
 For a Lie subalgebra $\h \subset \g$ the 
 subcomplex  of \textit{relative} Lie cochains, cf. Chapter 1.3 in \cite{Fu1}, $C^\bullet (\g, \h; V)$ give by 
 the condition that  both $\alpha$ and $d_{Lie}(\alpha)$ vanish when one of their
 arguments is in $\h$. Its cohomology groups are denoted by 
 $H^\bullet(\g, \h; V)$. 
 One source of relative Lie cocycles, see Section 2.2 of \cite{NT} arises from an $ad(\h)$-invariant
 projection $pr: \g \to \h$ and its curvature
 \begin{equation}
 \label{curvature}
 C(u \wedge v) = [pr(u), pr(v)] - pr ([u, v]): \Lambda^2 \g \to \h.
 \end{equation}
Assume for simiplicity that the $\g$ action on $V$ is trivial. 
Then for any $\h$-invariant polynomial $S \in Sym^l(\h^*)^{\h} \otimes V$ the cochain 
$\rho(S) \in C^{2l}(g; V)$ defined by 
$$
\rho(S)(v_1 \wedge \ldots \wedge v_{2l} )= \frac{1}{l!} \sum_{\sigma \in S_{2l}, \sigma(2i-1) < \sigma(2i)} (-1)^\sigma S(C(v_{\sigma(1)}, v_{\sigma(2)}), \ldots, 
C(v_{\sigma(2l-1)}, v_{\sigma(2l)}))
$$
is relative with respect to $\h$, closed, and its relative cohomology class is independent 
on the choice of the projection $pr: \g \to \h$. This defines the Chern-Weil homomorphism 
$$
\rho: Sym^\bullet(\h^*)^{\h} \to H^{2 \bullet}(\g, \h; k)
$$
We will need the following examples of relative cocycles:

\begin{enumerate}
	
\item When $\h \simeq \mathfrak{gl}_e (k)$, set $ch_{Lie} = \rho(tr(exp(x)))$ 
and $c_{1, Lie} = \rho(tr(x))$ for $x \in \h$;

\item When $h \simeq \mathfrak{sp}_{2p}(k)$, set $\widehat{A}_{Lie} = 
\rho\Bigg(\det \Big(\frac{y/2}{sinh(y/2)} \Big)^{1/2}\Bigg)$  for $y \in \h$. 

\item For a central extension of Lie algebras 
$$
0 \to \mathfrak{a} \to \widetilde{\g} \to \g \to 0
$$
and a $k$-vector space splitting $\widetilde{\g} \simeq \mathfrak{a} \oplus \g$,
the 2-cocycle $C: \Lambda^2 \g \to \mathfrak{a}$ is  the curvature as above.
In the cases we consider the 
cocycle may be chosen in $C^2(\g, \h; \mathfrak{a})$. 

\end{enumerate}

\subsection{Torsors over Harish-Chandra pairs and characteristic classes.}

\textbf{Definition.} A \textit{Harish-Chandra pair} $(\g, F)$ consists of a Lie algebra $\g$,
a (pro)algebraic group $F$ over $k$, an embedding of Lie algebras 
$\mathfrak{f} = Lie(F) \subset
\g$ and an action of $F$ on $\g$ which extends the adjoint action of $F$ on
$\mathfrak{f}$. 
A \textit{module} over a Harish-Chandra pair $(g, F)$ is an $F$-module $V$ with an 
$F$-equivariant Lie
morphism $\g \to End_k(V)$ extending the tangent Lie morphism on 
$\mathfrak{f}$. 

\medskip
\noindent
In this paper, $\mathfrak{f}$ will have
finite codimension in $\g$ and $F \simeq L \ltimes U$ with $L$ a finite dimensional reductive group and $U$
 a pro-unipotent infinite dimensional algebraic group. 
 
\bigskip
\noindent
\textbf{Definition.} A \textit{Harish-Chandra torsor} or a  flat $(\g, F)$-torsor over a scheme $Y$ is 
an $F$-torsor $\pi: P \to Y$ with an $F$-equivariant $\g$-valued 1-form $\gamma: T_P \to \g \otimes_k 
\mathcal{O}_P$ which restricts to the canonical Maurer-Cartan form (with values 
in $\mathfrak{f} \subset \g$) on the vector fields tangent to the fibers of $\pi: P \to Y$, and 
satisfies the Maurer-Cartan equation 
$$
d \gamma + \frac{1}{2} [\gamma, \gamma] = 0;
$$
where $d$ is the de Rham differential on $P$ and the bracket is computed in $\g$. In the 
infinite dimensional case some care must be taken to define such torsors. One possible
approach is to follow the pattern in \cite{Ye}, \cite{VdB} and work with representable 
functors. In our case, eventually we will only need direct images of differential forms 
from $P$ to $Y$, and since $F$ is a limit of affine groups, all geometric objects on $Y$ can
be defined using shaves on $Y$ with a coaction of functions on $F$, and so on. See 
Sections 2 and 3 in \cite{DGW} for a closely related case.

\medskip
\noindent 
Assume that $\gamma: T_P \to \g \otimes_k \mathcal{O}_P$ is onto and 
has kernel $T_\gamma$ of
finite constant rank $q$. Note that since $\gamma$ is injective on vertical vector fields, the
differential of $\pi: P \to X$ is injective on $T_{\gamma}$. Moreover, since 
$\gamma$ is $F$-equivariant, for points $x,y$ in the same orbit of $F$, the images
$(T_\gamma)_x$ and $(T_\gamma)_y$ in $(T_Y)_{\pi(x) = \pi(y)}$ agree. Let 
$\mathcal{F}  \subset T_Y$ the resulting rank $q$ sub-bundle   on $Y$. 
\begin{lemma}
	In the situation described, let $T_\pi = Ker(d\pi) \subset T_P$ be the vector fields 
	tangent to the fibers. Then the sub-bundles $\mathcal{F} \subset T_Y$ 
	and $\mathcal{G} = T_\pi \oplus ker(\gamma) \subset T_P$ are integrable (i.e. stable under the Lie bracket
	of vector fields). 
	Let $\Omega_{\mathcal{F}} \subset \Omega^1_Y $ be the annihilator of $\mathcal{F}$ and 
	$\Omega_{\mathcal{G}} \subset \Omega^1_P$ the annihilator of $\mathcal{G} $. 
	Denote by $F^r\Omega^\bullet_Y$, resp. $F^r\Omega^\bullet_P$, the graded ideal generated
	by the $r$-th power of $\Omega_{\mathcal{F}}$, resp. the $r$-th power of $\Omega_{\mathcal{G}}$ .
	Then both $F^r\Omega^\bullet_Y$ and $F^r\Omega^\bullet_P$ are preserved by the corresponding
	de Rham differentials 
	and there is a morphism of complexes of sheaves $F^r \Omega^\bullet_Y 
		\to \pi_* F^r \Omega^\bullet_P$
\end{lemma} 
 \textit{Proof:}  We start with $\mathcal{G}$. 
 If $v_1$, $v_2$ are two vector fields in $Ker(\gamma)$ then formula (\ref{differential})
 for $d\gamma$ and the Maurer-Cartan equation for $\gamma$ imply that the bracket 
 $[v_1, v_2]$ is also annihilated by $\gamma$.  The fact that vector fields tangent to the fibers
 are closed with respect to the Lie bracket holds for any smooth $\pi$. Finally, let's assume
 that $v_1 \in Ker(\gamma)$ and $v_2 \in T_\pi$.  Then the quadratic term in the Maurer-Cartan 
 equation vanishes on $v_1 \wedge v_2$ (as $\gamma(v_1) = 0$) and we are left with 
 $$
 0 = d \gamma (v_1 \wedge v_2) = v_1 \cdot \gamma(v_2) - v_2 \gamma(v_1)  - \gamma ([v_1, v_2]). 
 $$
 The second term is zero by assumption on $v_1$ and the first is in $\mathfrak{f} \otimes_k \mathcal{O}_P
 \subset \g \otimes_k \mathcal{O}_P$. Hence the third term is in $\mathfrak{f} \otimes_k \mathcal{O}_P$ as
 well, which means $[v_1, v_2] \in \gamma^{-1} (\mathfrak{f} \otimes_k \mathcal{O}_P) = \mathcal{G}$.

 In particular, a bracket of two  $F$-equivariant vector fields in $Ker(\gamma)$ is 
 again an $F$-equivariant vector field in $Ker(\gamma)$. Its $F$-equivariant descent is
  a rank $q$ subbundle in the Atiyah algebra $At_P$
 of $P$ on $Y$
 (= the $F$-equivariant descent of all vector fields on $P$), which is also closed under  Lie bracket. 
 It projects isomorphically to a sub-bundle $\mathcal{F} \subset T_Y$ (as 
 $Ker(\gamma)$ has trivial intersection with $T_\pi$) which is closed with respect to the Lie bracket since
 $At_P \to T_P$ is compatible with brackets. 
 By construction, 1-forms on $Y$ which vanish on $\mathcal{F}$ pull back to $F$-equivariant 
 1-forms on $P$
 which vanish on $\mathcal{G}$. 
 Hence the morphism of sheaves of dg algebras $\Omega^\bullet_Y \to \pi_* \Omega^\bullet_P$
 is compatible with multiplicative filtrations $F^r \Omega^\bullet$ induced by the two foliations. 
  $\square$
 
\bigskip
\noindent 
\textbf{Definition.}  In the situation of the previous lemma we will say that the Harish-Chandar torsor
$(P, \gamma)$ is \textit{foliated} over $\mathcal{F} \subset T_Y$. We will see that in this case some of
its characteristic classes what apriori belong to $H^\bullet_{DR}(Y)$  admit a lift to cohomology 
of $F^r \Omega^\bullet_Y$.

\bigskip
\noindent 
Let $V$ be a Harish Chandra module over $(\g, F)$ with trivial action (as will be in our applications). 
For a Lie $l$-cochain $\alpha: \Lambda^l \g \to V$ the composition $\Lambda^l T_P
\to \Lambda^l \g \otimes_k \mathcal{O}_P \to V \otimes_k \mathcal{O}_P$ may be viewed
as a $V$-valued $l$-form on $P$ and the Maurer-Cartan equation ensures that the 
resulting \textit{Gelfand-Fuks morphism} cf. Chapter 3.1.C,D in \cite{Fu1} 
also agrees with differentials:
\begin{equation}
	\label{gf}
GF: C^\bullet(\g; V) \to \Gamma(P, \Omega^\bullet_P \otimes_k V).
\end{equation} 
We want to use this observation to study characteristic
classes of $P$ in the cohomology of $Y$. 
Note that to obtain classes on $Y$ we need to work with objects which are invariant with 
respect to the reductive subgroup of $F$.

For a subalgebra $\h \subset \mathfrak{f} = Lie(F) \subset \g$, relative cochains in $C^\bullet(\g, \h; V)$ map to 
the subcomplex $(\Omega^\bullet \otimes_k V)_{\h-basic}$  of $\h$-basic forms, 
i.e. forms $\beta$ such 
that $L_v \beta = 0$ and $\iota_v \beta = 0 $ for any $v \subset \h$ 
(we use the same letter $v$ for 
the vertical vector field on $P$ induced by $v$ via the action of $F$).

\begin{lemma} 
	\label{filtrations} 
	Assume that $F = U \rtimes  H$ is  semi-direct product of a finite dimensional connected reductive group $H$ and a
	pro-unipotent group $U$. If $\h = Lie(H)$, there exists a quasi-isomorphism of sheaves of dg-algebras
	$$
	\Omega^\bullet_Y \to \pi_* \Big( (\Omega^\bullet_P)_{\h-basic}  \Big)
	$$
	If $P$ is foliated over $\mathcal{F} \subset T_Y$ then for any $r \geq 0$ the natural morphism of ideal sheaves
	$$
	F^r \Omega^\bullet_Y \to \pi_* \Big( (F^r \Omega^\bullet_P)_{\h-basic} \Big) 
	$$
	is a quasi-isomorphism of complexes of sheaves. 
\end{lemma}
\textit{Sketch of Proof.} First assume that $U$ is finite dimensional and look at the first 
quasi-isomorphsm. Since $H$ is connected, the pushforward 
of $\h$-basic forms from $P$ to $Y$ may be identified with the pushforward of forms on $P/H$ to $Y$.
But $P/H \to Y$ is a bundle with affine fibers so the assertion follows from the relative
Poincare Lemma
(triviality of relative de Rham cohomology for fibrations by affine spaces). For 
infinite dimensional pro-unipotent $U$, we first consider the finite dimensional unipotent factors 
and then pass to a limit, as in Theorem 6.7.1 in 
\cite{VdB} .

In view of the unfiltered 
quasi-isomorphism, its filtered version reduces to showing that
the maps induced on associated graded quotients are quasi-isomorphisms. First, 
$$
F^r \Omega^\bullet_Y /
F^{r+1} \Omega^\bullet_Y \simeq \Lambda^r(Q) \otimes \Lambda^\bullet \mathcal{F}^\vee
$$
 where 
$Q = T_Y/\mathcal{F}$ and we consider $\Lambda^\bullet \mathcal{F}^\vee$ as a complex of sheaves
with the differential similar to that of in formula
(\ref{differential2}) (in other words, it is the de Rham differential 
of the Lie algebroid $\mathcal{F} \subset T_Y$).  On the other hand, since the pullback of $\mathcal{F}$
is isomorphic to $Ker(\gamma) \subset \mathcal{G}$ and the pullback of $Q$ is isomorphic to $T_P/\mathcal{G}$, 
by the projection formula we get
$$
\pi_* \Big( (F^r \Omega^\bullet_P)_{\h-basic} \Big) / \pi_* \Big( (F^{r+1} \Omega^\bullet_P)_{\h-basic} \Big) 
\simeq \Lambda^r(Q) \otimes \Lambda^\bullet \mathcal{F}^\vee
\otimes_{\mathcal{O}_Y} \pi_* \Big( (\Omega^\bullet_\pi)_{\h-basic} \Big) 
$$
where $\Omega^\bullet_\pi$ is the relative de Rham complex of $\pi: P \to Y$. Hence we just
need to show that
$$
\mathcal{O}_Y \to \pi_* \Big( (\Omega^\bullet_\pi)_{\h-basic} \Big) 
$$
is a quasi-isomorphism, i.e. which again follows from
 the relative Poincare Lemma. $\square$
 
 \bigskip
 \noindent
 In view of the previous result, passing to cohomology in (\ref{gf}) we obtain 
 an $\h$-relative version of Gelfand-Fuks map, which we denote also by $GF$:
 $$
 GF: H^\bullet_{Lie}(\g, \h; V) \to H^\bullet_{DR}(Y) \otimes_k V.
 $$
 Recall that whenever we use $GF$ we assume that $\g$-action on $V$ is trivial, otherwise
 the right hand side would involve the de Rham cohomology of associate vector
 bundle $V_P$ with a flat connection induced by the Harish-Chandra module structure.

 \begin{proposition}
 	\label{compare-classes} 
 	For $V = k$ in the setting of Lemma \ref{filtrations}  
 	the composition of the Lie algebraic Chern-Weil map and the 
 	Gelfand-Fuks map
 	$$
 	Sym^\bullet(\h^*)^\h \to C^{2\bullet}(\g, \h; k) \to H^{2\bullet}_{DR} (Y)
 	$$
 	is the classical Chern-Weil map of the torsor $P_H$, associated to $P$ via the group 
 	homomorphism $F \to H$. If $P$ is foliated over $\mathcal{F} \subset
 	T_Y$, for every $r \geq 0$  the composition admits a canonical lift
 	$$
 	Sym^r(\h^*)^\h \to H^{2 r}_{DR} (Y, F^r \Omega^\bullet_Y)
 	$$
 	\end{proposition}
 \textit{Proof.} Let $P_U = P/H$ then $P \to P_U$ maybe be viewed as the $H$-torsor
 pulled back from $Y$ via $P_U \to Y$. Since $P_U \to Y$ induces isomorphism on 
 de Rham cohomology (relative Poincare Lemma)  we can replace $Y$ by $P_U$ and
 assume that $U$ is trivial. Then the composition 
  $$
  \nabla = pr \circ \gamma: T_P \to \g \otimes_k \mathcal{O}_P \to \h \otimes_k \mathcal{O}_P
  $$
  is a connection on $P \to P_U$ and the Lie theoretic curvature $C: \Lambda^2 \g 
  \to \h$ gives the classical curvature $R_\nabla = \Lambda^2 \gamma
  \circ C: \Lambda^2 T_P \to \h \otimes_k 
  \mathcal{O}_P$. The assertion follows since the classical Chern-Weil map may be
  computed by evaluating invariant polynomials on $R_\nabla $. 
  
  For the second assertion, we observe that relative Lie cochains define a global section of
  $\pi_* \Big( (F^r \Omega^\bullet_P)_{\h-basic} \Big) $ and the result follows by 
  application of Lemma \ref{filtrations}. 
  $\square$

\subsection{Algebraic version of the Index Theorem.}

Consider the formal Weyl algebra 
$\mathcal{D}_p$, the completion (at the augmentation ideal) 
of the universal enveloping 
of the Heisenberg Lie algebra with generators $x_1, \ldots, x_p, y_1, 
\ldots, y_p, h$
and the only nontrivial commutators given by $[y_j, x_i] = \delta_{ij}h$. 
In other words, as a vector space $\mathcal{D}_p$ is isomorphic to 
$k[[x_1, \ldots, x_p y_1, \ldots, y_p, h]]$ but with 
nontrivial commutation relations $y_i x_i = x_i y_i + h$. 

We consider the associative algebra $\mathcal{E} = \mathfrak{gl}_e(\mathcal{D}_p)$ and the Lie
algebra $Der(\mathcal{E})$ of its  continuous $k[[h]]$-linear derivations. Of course, some of the derivations are
are inner and hence there is a Lie morphism $\mathcal{E} \to Der(\mathcal{E})$. In addition, 
any commutator in $\mathcal{D}_p$ is divisible by $h$, so commuting with 
a scalar $\frac{1}{h} \mathcal{D}_p$-valued matrix also gives a derivation 
of $\mathcal{D}_p$. We claim a short exact sequence
\begin{equation}
	\label{dersofE}
0 \to \frac{1}{h} k[[h]] \to \Big(\frac{1}{h} \mathcal{D}_p + \mathcal{E}\Big)
 \to  Der(\mathcal{E}) \to 0
\end{equation} 
For $e =1$ this is well-known, see e.g. Section 2.3 in \cite{BNT} and equation (3.2) in \cite{BK}.
For general $e$ follows from the fact that the
quotient of all derivations by inner derivations (i.e. first Hochschild cohomology group) is a
Morita invariant, cf. Chapter 1.2 in \cite{Lo}. 

The Lie algebra $\g = \Big(\frac{1}{h} \mathcal{D}_p + \mathcal{E} 
\Big)$ has a reductive subalgebra $\mathfrak{gl}_e \oplus \mathfrak{sp}_{2n}$ 
(matrices with values in $k \subset \mathcal{D}_p$ plus an isomorphic copy of 
$\mathfrak{sp}_{2n}$ in 
$\frac{1}{h}\mathcal{D}_p$ spanned by commutators  of $1, \frac{1}{h} x_i x_j, \frac{1}{h} y_i x_j, 
\frac{1}{h} y_i y_j$).  We also
consider the abelian subalgebra $\mathfrak{a} = \frac{1}{h} k[[h]]$.  To take into 
 account  $\mathfrak{a}
\cap \mathfrak{gl}_e = k$, introduce $\mathfrak{a}' \subset \mathfrak{a}$
with topological basis given by $h^{i}$, $i = -1, 1, 2, \ldots$. Then 
$$
\h = \mathfrak{gl}_e \oplus \mathfrak{sp}_{2n} \oplus \mathfrak{a}'
$$
is a Lie subalgebra of $\g$. By Hochschild-Serre spectral sequence
$$
H^\bullet(\g, \h; V) \simeq H^\bullet( Der(\mathcal{E}), 
(\mathfrak{gl}_e/k )\oplus \mathfrak{sp}_{2p}; V)
$$
where $V$ is a module over $Der(\mathcal{E})$  (and thus also 
a module over $\g$). 
Below, we are interested in the cohomology of $Der(\mathcal{E})$
but find it easier to do computations in  $\g$.

\bigskip
\noindent
Now we would like to state a vanishing 
lemma for homogeneous components of a
particular class in $H^\bullet(\g, \h; k((h)))$. Its proof will take the rest of the section. 
A reader willing
to treat it as a black box may wish skip to Section 4. We follow the notation and 
exposition in \cite{GLX} which deals with a version of Algebraic Index Theorem 
that is most convenient for our setting. 

This class can be defined via the Chern-Weil construction.  
To fix a projection $pr: \g \to \h$  we
introduce a filtration on $\g$ 
by giving $h$ degree 2 and $x_i, y_j$ degree 1. Since elements of $\g$ involve
infinite sums, $\g$ splits into a direct product
$\prod_{i \geq -2} \g_i$.  For the first two factors in $\h$, we project $\g$
onto $\g_0 \simeq \mathfrak{gl}_e \oplus \mathfrak{sp}_{2n}$ (recall that 
$\mathfrak{sp}_{2n}$ is spanned by the commutators in the degree zero part
$\big(\frac{1}{h}\mathcal{D}_p\big)_0$)
and set $k = \frac{h}{h} k$ to be in the kernel on the projection onto $\mathfrak{sp}_{2n}$.
For $\mathfrak{a}'$ we choose any projection $\frac{1}{h} \mathcal{D}_{p} \to \mathfrak{a}'$
and extend it by 0 to trace zero matrices in $\mathcal{E}$.

\begin{lemma}
	\label{vanishing}
	For $\h = \mathfrak{gl}_e \oplus \mathfrak{sp}_{2p} \oplus \mathfrak{a}'$, let  $ch_{Lie}(\mathfrak{gl}_e) $,
	  $\widehat{A}_{Lie}(\mathfrak{sp}_{2p})$ and $C(\mathfrak{a}')$
	  be the classes induced by the Chern-Weil constrution at the end of Section 3.1, from the respective factors.  Let
	$$
	 \tau_{\mathcal{D}_p} = ch_{Lie}(\mathfrak{gl}_e) \widehat{A}_{Lie}(\mathfrak{sp_{2n}})  exp(- C(\mathfrak{a}')) \in 
	 H^{even}_{Lie} (\g, \h; k((h)))
	$$
	Then the components of $\tau_{\mathcal{D}_p}$ of degree $> 2p$ are equal to zero. 
	
\end{lemma}
\textit{Proof.}  Our proof is based on the fact that the class of the theorem arises from 
the study of periodic cyclic homology of the associative 
algebra $A = gl_e(\mathcal{D}_p)
\otimes_{k[[h]]} k((h))$. The vanishing  will
follow from the fact that this class lifts to the negative cyclic homology of $A$. 
We briefly recall the relevant definitions here, omitting details that do not contribute to the 
proof.

\bigskip
\noindent 
If $A$ is a unital $k$-algebra, set $\overline{A} = A/k \cdot 1$ and $C_{-l} (A) = A 
\otimes \overline{A}^{\otimes l}$. This is the cohomological grading, rather than homological grading
used in some sources, although the indices are written as subscripts to avoid confusion with the
Hochschild cohomology complex. 
The standard formulas e.g. in 1.1.1 and 2.1.8 of \cite{Lo} define the Hochschild and cyclic differentials,  
$$
b: C_{\bullet}(A) \to C_{\bullet+1} (A); \qquad B: C_{\bullet} (A) \to C_{\bullet -1} (A)
$$
that satisfy $b^2 = 0, B^2 = 0, Bb + bB = 0$. 
Introducing a formal variable $u$ of cohomological degree 2,
consider two complexes
$$
CC_\bullet^-(A) = (C_\bullet[[u]], b + uB), \qquad CC_\bullet^{per} = (C_\bullet((u)), b + u B)
$$
with cohomology defining the negative cyclic homology $HC^{-}_\bullet(A)$ and the periodic cyclic
homology $HC^{per}_\bullet(A)$, 
respectively. 
In both cases $1 \in A = C_0(A)$ satisfies $(b + u B)(1) = 0$ and thus gives a cohomology class.

\bigskip
\noindent
In the case $A = gl_r(\mathcal{D}_p) \otimes_{k[[h]]} k((h))$, by the work of Shoikhet, cf. 
\cite{Sh} and Section 1.2 in \cite{CFW}, there is  a quasi-isomorphism 
$(C_\bullet(A), b) \to (\Omega^{-\bullet}((h)), h L_\pi)$
where
$\Omega^{- \bullet}$ stands for formal differential forms in $x_1, \ldots, x_p, 
y_1, \ldots, y_p$ and $\pi = \sum (\partial/\partial x_i) \wedge (\partial/\partial y_i)$ is 
the standard Poisson bivector.
Furthermore, Willwacher proved in \cite{Wi}, that this quasi-isomorphism sends $B$ to the de 
Rham differential. This results in a quasi-isomorphism 
$$
(C_{\bullet}(A)((u)), b+uB) \simeq (\Omega^{-\bullet}((h))((u)), h L_\pi + u d_{DR}).
$$
On the left hand side, $L_{\pi} = d \iota_{\pi} + \iota_{\pi} d: \Omega^{l} \to \Omega^{l-1}$ 
. 

\bigskip
\noindent 
The explicit construction on the above quasi-isomorphism is quite non-trivial and involves
integration over configuration spaces of points. The Lie  algebra $Der(\mathcal{E})$ 
of derivations of $\mathcal{E} = \mathfrak{gl}_e(\mathcal{D}_p)$ 
acts on periodic and negative cyclic complexes. It follows from the explicit construction in \cite{Sh} 
 that the quasi-isomorphsm 
 is compatible with the action of $\mathfrak{pgl}_e(k) \oplus \mathfrak{sp}_{2n} \subset
Der(\mathcal{E})$ but
not the full algebra of derivations.
  However, it
can be upgraded to a cocycle in the relative Lie algebra cohomology of $Der(\mathcal{E})$. 
As before, we work
with $\g$ (a central extension of $Der(\mathcal{E})$) rather than derivations of $A$. Then 
the upgraded quasi-isomorphism is an element of total degree 0
$$
\tau_{Lie} \in C^\bullet (\g, \h; Hom_{k((h))} (( CC_{\bullet}^-(A), b + uB), (\widehat{\Omega}^{-\bullet}((h))[[u]], 
h L_\pi + u d_{_{DR}}))
$$
We note here that the action of $\g$ is not $\widehat{\Omega}^{-\bullet}((h))$ is not $h$ linear
(in fact the element $h \in \g$ acts by zero, so the action factors through the quotient by $h$), and
all cochains are only $k$-linear maps. 
This is usually considered for the periodic cyclic complex but we emphasize that
at this step no inversion of $u$ is necessary. 
Next, one can construct two $\g$-invariant homomorphisms 
$$g_u \circ e^{-\frac{h}{u} \iota_{\pi}} , g_h \circ e^{- \frac{u}{h} \omega}  \in 
 Hom_{k((h))} ((\widehat{\Omega}^{-\bullet}((h))[[u]], 
h L_\pi + u d_{_{DR}}), (\widehat{\Omega}^{\bullet}((h))((u)), d_{_{DR}}))
$$
which are homotopic in the full Lie algebra cohomology complex. 

The first homomorphism needs inversion of  $u$. We first use the identity 
 $e^{\frac{h}{u} \iota_{\pi}} (h L_\pi + u d_{_{DR}})
= (u d_{_{DR}})e^{\frac{h}{u} \iota_{\pi}} $ to land in the complex 
$(\widehat{\Omega}^{-\bullet}((h))((u)), ud_{_{DR}})$. Then we
apply the re-grading operator $g_u$ which is
an isomorphism of complexes
$$
(\widehat{\Omega}^{-\bullet} ((h))((u)), u d_{_{DR}} ) \to 
(\widehat{\Omega}^\bullet((h))((u)), d_{_{DR}}) 
$$
sending an $i$-form $\alpha$ to $u^{-i} \alpha$. Thus, on the left hand side $\alpha$ has cohomological degree
$(-i)$, on the right hand side degree $i$ and adjustment by $u^{-i}$ makes the re-grading a degree 0 operator. 

The second 0-cocycle  uses $h^{-1}$ and the formal symplectic form
$\omega = \sum dx_i \wedge dy_i$. In this case, we start with the identity 
$e^{- \frac{u}{h} \omega}(h L_\pi + u d_{_{DR}}) = h L_\pi e^{- \frac{u}{h} \omega}$
to land in the complex $(\widehat{\Omega}^{-\bullet}((h))[u^{-1}, u]], h L_\pi)$ 
and then apply a re-grading isomorphism $g_h$
$$
(\widehat{\Omega}^{-\bullet} ((h))((u)), h L_\pi ) \to 
(\widehat{\Omega}^\bullet((h))((u)), d_{_{DR}} ) 
$$
which sends  an $i$-form $\alpha$ to $\frac{h^{i-n}}{u^n} *_\omega (\alpha)$. Here, 
the symplectic Hodge operator $*_\omega: \Omega^i \to \Omega^{2n-i}$ is defined
by $\beta \wedge *_\omega (\alpha) = \omega^n \langle \beta, \alpha \rangle_\pi$ 
and $\langle \cdot, \cdot \rangle_\pi$ is
the pairing on $i$ forms induced by $\pi$. 

\bigskip
\noindent
With these preliminaries, we now prove the vanishing claimed in Lemma \ref{vanishing} in 
several steps

\bigskip
\noindent
\textit{Step 1.} First consider the class 
	$$
	\tau_1 = \sum \tau_1^{k, l} u^l \in \bigoplus C^{2k} (\g, \mathfrak{h}; \widehat{\Omega}^{-2(k+l)}((h)) u^l) 
	$$ 
	is obtained by pairing $1 \in C^0_{Lie}(\g, \mathfrak{h}; 
	(CC^-_0(A), b + uB))$ with $\tau_{Lie}$ (we recall that the grading on differential forms has been inverted 
	 at this step). The range of indices is $k \geq 0$ and $ 0 \leq k + l \leq n$. 
	Then move on to the class $\tau_2 = e^{-\frac{h}{u} \iota_{\pi}}\tau_1$ which lives in the same complex but
	with the differential $ud_{_{DR}}$ instead of $h L_{\pi} +u d_{_{DR}}$. The components $\tau_2^{k, l} u^l$ of the class
	$\tau_2$ can be nonzero in the same range $0 \leq k, 0 \leq k+l \leq n$. 
	
	\bigskip
	\noindent
	\textit{Step 2.} Now consider the regraded class 
	$$
	\tau_3 = g_u \tau_2 = \sum \tau_2^{k, l} u^{-2k -2l} u^l \in \bigoplus C^{2k} (\g, \mathfrak{h}; \widehat{\Omega}^{2(k+l)}((h)) u^{-2k-l}) 
	$$ 
	(where the differential forms now have the usual grading)
	and observe that
	the exponents of $u$ are in the range $(-2n, \ldots, 0)$.  The next step is to replace the de Rham complex
	$(\widehat{\Omega}^\bullet, d_{_{DR}})$  by a quasi-isomorphic complex $(k, 0)$. Note that the 
	projection $\widehat{\Omega}^\bullet \to k$ (which vanishes on forms of positive degrees and sends a power
	series in degree zero to its constant term) - is not $\g$-equivariant. However, it can be extended to 
	a quasi-isomorphism of complexes
	$\widetilde{f}: C^\bullet(\g, \mathfrak{h}; \widehat{\Omega}^\bullet) \to 
	C^\bullet(\g, \mathfrak{h}; k) $ using 
	Lemma \ref{homotopies} (i) below. This leads to a class 
	$$
	\tau_4 = \sum \tau^{4k+2l}_4 u^{-(2k+l)} = \widetilde{f}(\tau_3) \in \bigoplus C^{4k+2l} (\g, \mathfrak{h}; k((h)) u^{-2k-l}) 
	$$
	Note that $0 \leq k, 0 \leq k+l \leq n$  imply $0 \leq 4k + 2l \leq 4n$. This agrees with the 
	fact that  the de Rham cohomology of $M$ is nonzero in the range $(0, \ldots, 4n)$. 
	Now we use the Algebraic Index Theorem in the form proved in 
	\cite{GLX} to claim that $\tau_4 =\tau_{\mathcal{D}_p}$, the class introduced
		in the lemma. 
		
\bigskip
\noindent
\textit{Step 3.}  We are finally in position to prove that on the level of cohomology the coefficients $\tau_4^{4k + 2l}$
	vanish when $4k + 2l > 2n$. Indeed, instead of 
	$$
	\tau_3 = g_u e^{- \frac{h}{u} \iota_\pi} \tau_1 
	$$
	we could consider 
    $$
    \tau_3' = g_h e^{- \frac{u}{h} \omega} \tau_1
    $$
    which has trivial coefficients of $u^{-(2k+l)}$, $2k+l \geq n$ since the only negative power of 
    $u$ created is the factor $u^{-n}$ in the definition of $g_h$.

 \bigskip
 \noindent
 \textit{Step 4.}   
    To show that $\tau_3$ and $\tau_3'$ have the same cohomology class we use the identity 
    $$
    e^{\frac{\omega}{uh}} \circ (g_h \circ e^{- \frac{u}{h} \omega})
    = g_u \circ e^{\frac{h}{u} \iota_{\pi}};
    $$
    which can either be established by direct computation, or by using a basic observation of Hodge Theory 
    that operators $\wedge \omega $ and $\iota_\pi$ generate an $\mathfrak{sl}_2$ action on the de Rham 
    complex, hence the above identity can be obtained as the image of the group
    level identity in $SL_2$ 
    $$
\begin{pmatrix}
	1 & 0 \\
	h & 1
	\end{pmatrix} \cdot
\begin{pmatrix}
	0 & h^{-1} \\
	-h & 0
\end{pmatrix} \cdot
\begin{pmatrix}
	1 & 0 \\
	h & 1
\end{pmatrix} =
\begin{pmatrix}
	1 & h^{-1} \\
	0 & 1
\end{pmatrix} 
    $$
    under the action homomorphism.
    
  Finally, we note that
     $e^{\frac{\omega}{uh}} $ is homotopic to identity. This follows from the fact that 
    $\omega  = d_{_{DR}} (\alpha)$ where $\alpha = \frac{1}{2} \sum (x_i d y_i - y_i dx_i)$ is
    the Euler vector field converted to 1-forms using $\omega$. Therefore we have 
    $d_{_{DR}} \varphi + \varphi d_{_{DR}} = e^{\frac{\omega}{uh}}$ where 
    $$
    \varphi = \sum_{k \geq 0} \frac{\alpha}{uh} \frac{\omega^k}{(uh)^k (k+1)!}.
    $$
    Again, the homotopy does not agree with the $\g$-action (only with the $\mathfrak{h}$-action) hence
    we use  Lemma \ref{homotopies}(ii) below to obtain a homotopy 
    between identity operator and $e^{\frac{\omega}{uh}}$.
    $$
    \widetilde{\varphi}: C^\bullet(\g, \mathfrak{h}; \widehat{\Omega}^\bullet((h)) ((u)))
    \to C^{\bullet-1} (\g, \mathfrak{h}; \widehat{\Omega}^\bullet((h)) ((u)))
    $$
    Hence we can 
    use the above class $\tau_3'$ instead of $\tau_3$. Since $\tau_3'$ by construction has 
    at worst poles of order $\leq n$ in the $u$ variable, the assertion follows.  $\square$
    
\bigskip
\noindent     
We finish here with a homotopy lemma used  above
\begin{lemma}
	\label{homotopies}
Let $M^\bullet$ and $N^\bullet$ two complexes of modules over a Lie algebra $\g$ and 
$f: M^\bullet \to N^\bullet$, $g: N^\bullet \to M^\bullet$  chain maps
compatible with an action of a subalgebra 
$\h \subset \g$.  	Suppose that $\varphi: M^\bullet \to M^{\bullet-1}$ is 
a homotopy satisfying $d \varphi + \varphi d = gf - 1_M$. Denote by 
$$
\delta: C^\bullet(\g, \mathfrak{h}; \cdot) \to C^\bullet(\g, \mathfrak{h}; \cdot) \textrm{ and }
d_{Hom}: C^\bullet(\g, \mathfrak{h}; \cdot) \to C^\bullet(\g, \mathfrak{h}; \cdot)
$$
the first terms in the formula (\ref{differential2}) and the combination of the second 
term of (\ref{differential2}) with the internal differential on $M^\bullet, N^\bullet$, respectively. 
Let 
$$
f_{Hom}: C^\bullet(\g, \mathfrak{h}; M^\bullet) \to C^\bullet(\g, \mathfrak{h}; N^\bullet),
g_{Hom}: C^\bullet(\g, \mathfrak{h}; N^\bullet) \to C^\bullet(\g, \mathfrak{h}; M^\bullet)
$$
be the morphisms of complexes with $d_{Hom}$ differentials, induced by $f, g$, respectively, 
and $\varphi_{Hom}: C^\bullet(\g, \mathfrak{h}; M^\bullet)  \to C^\bullet(\g, \mathfrak{h}; M^{\bullet-1})$ a homotopy induced by $\varphi$.  

(i) If $g$ is compatible with $\g$ action and the 
side conditions $\varphi \varphi = 0$, $\varphi g = 0$, $f \varphi = 0$ hold, then 
$$
\widetilde{f} = f_{Hom} (1 + (\delta \varphi_{Hom}) + (\delta \varphi_{Hom})^2 + (\delta \varphi_{Hom})^3 + \ldots) 
$$
is compatible with Lie algebra cohomology differentials $d_{Lie} = d_{Hom} + \delta$. 

(ii) If $M^\bullet = N^\bullet$, $g = 1_M$ and $f$ is compatible with the $\g$-action  then 
$f_{Hom}$ is compatible with $d_{Lie} = d_{Hom} + \delta$ and 
$$
\widetilde{\varphi} =\varphi_{Hom}(1 + (\delta \varphi_{Hom}) + (\delta \varphi_{Hom})^2 + (\delta \varphi_{Hom})^3 + \ldots) 
$$ 
is a homotopy between $1$ and $f_{Hom}$. 
	\end{lemma}
\textit{Proof.} 
Part (i) is a consequence of the Basic Perturbation Lemma, cf. \cite{Ma}
Part (ii) is easier to establish by direct computation although it is also a very degenerate case of
the Ideal Perturbation Lemma, cf. \cite{Ma}. $\square$

\section{Proof of the main result.}

In this section we prove  Theorem \ref{main} by studying the characteristic 
class 
$$
\tau_Y(E) = exp(- \frac{c_1(N)}{2}) ch(E)  \widehat{A}(B)   exp(-c(\mathcal{O}_h)|_Y) \in H^\bullet_{DR}(Y)((h))
$$
For $y \in Y$ the preimage of the maximal ideal $\mathfrak{m}_y \subset \mathcal{O}_{Y, y}$ in the 
stalk at $y$, with respect to the reduction mod $h$ map $\mathcal{O}_{h, y} \to \mathcal{O}_{Y, y}$, 
is a maximal ideal $\mathfrak{m}_{h, y} \subset \mathcal{O}_{h, y}$.
Adapting the classical proof of the Darboux theorem, we show that after completion at this maximal ideal the triple
$(\mathcal{O}_h, E_h, End_{\mathcal{O}_h}(E_h))$ is isomorphic - non-canonically! - to a 
similar triple independent of $y$ or $Y$. Different choices of isomorphisms will give the Harish Chandra
torsor $P_{\mathcal{D}, \mathcal{M}}$ inducing $\tau_Y(E)$ via the Gelfand-Fuks map. 
Further, it is actually lifted from a quotient torsor $P_{\mathcal{E}}$ 
and Theorem \ref{main} will be reduced to the study of the class $\tau_{\mathcal{D}_p}$ where
one uses  Lemma \ref{filtrations} and Lemma \ref{vanishing}. 

\subsection{Standard formal models: $\mathcal{D}, \mathcal{M}$ and $\mathcal{E}$.}

\bigskip
\noindent
Below
for $n = p+q$ we will assume that $\mathcal{D}_q$ is the Weyl algebra built on the variables $x_i, y_i, h$,
with $i = 1, \ldots, q$, that $\mathcal{D}_p$ corresponds to the values
$i = q+1, \ldots, q+p$ while $\mathcal{D}$ is
the Weyl algebra on the full set of variables with $i = 1, \ldots, p+q = n$. 
Fixing decomposition $n = p+ q$ and an integer $e \geq 1$, define 
a left $\mathcal{D}$-module 
$$
\mathcal{M}: =  \big[\mathcal{D}/\mathcal{D}\langle y_1, 
\ldots, y_q \rangle \big]^{\oplus e} \simeq \big(\mathcal{D}_q
/\mathcal{D}_q \langle y_1, \ldots, y_q \rangle \big) \widehat{\otimes}_{k[[h]]}
\mathcal{D}_{p}^{\oplus e}
$$
(on the right hand side, we use completed tensor product). 
The second presentation implies the following isomorphism for endomorphisms
of $\mathcal{M}$ (which we assume to be acting \textit{on the right}):
$$
\mathcal{E} := End_{\mathcal{D}} (\mathcal{M}) \simeq gl_e(\mathcal{D}_p)
$$

\begin{lemma} 
	\label{isoms}
	For any $Y, \mathcal{O}_h, E_h$ as before and $y \in Y$, denote by 
	$\widehat{(\ldots)}$ the completion of a stalk at $y$ with respect to
	the maximal ideal $\mathfrak{m}_{h, y}$. 
	
	(1) There exist compatible isomorphisms
	$$
	\sigma_{\mathcal{D}}: \widehat{\mathcal{O}_{h, y}} \to \mathcal{D}, 
	\qquad
	\sigma_{\mathcal{E}}:  \widehat{End}_{\mathcal{O}_{h, y}}(E_{h, y})\to \mathcal{E}
	\qquad 
	\sigma_{\mathcal{M}}: \widehat{E_{h, y}} \to \mathcal{M}
	$$
	as filtered algebras and modules, respectively. Moreover, 
	if an algebra isomorphism 
	$\sigma_\mathcal{D}$ admits a compatible module isomorphism $\sigma_{\mathcal{M}}$ 
	then $\sigma_{\mathcal{D}}$
	sends the completion $\widehat{I_h}$
	of $I_h = Ker(\mathcal{O}_h \to \mathcal{O}_Y)$
	to the double sided ideal $\mathcal{J} \subset
	\mathcal{D}$ generated by $y_1, \ldots, y_q, h$. 
	
	(2) If $\sigma_{\mathcal{D}}(\widehat{I_h}) = \mathcal{J}$ then 
	$\sigma_{\mathcal{D}}$ extends to a pair of compatible isomorphisms
	$(\sigma_{\mathcal{D}}, \sigma_{\mathcal{M}})$. 
	\end{lemma}
\textit{Proof.} 
Modulo $h$, we can construct an isomorphism of $\widehat{\mathcal{O}}_{M, y}$ and 
$k[[x_1, \ldots, x_n, y_1, \ldots y_n]]$ since $y$ is a smooth point, $k$ has characteristic
zero and $X$ has dimension 
$2n$. Since $Y$ is smooth, its ideal $I_Y \subset \mathcal{O}_{M, y}$ is generated by a regular
sequence and we can assume that the isomorphism sends the completion on $I_Y$ to the ideal
generated by $y_1, \ldots, y_q$. 

We can also adjust the isomorphism to be compatible with 
the symplectic forms. Let $\alpha$ be a two form on $Spec(k[[x_1, \ldots, x_n, y_1, \ldots, y_n]])$ induced
from $\omega$ via the initial isomorphism. Decompose it into  homogeneous components:
$\alpha = \alpha_0 + \alpha_1 + \alpha_2 + \ldots$ where each $\alpha_i$ is a two form with 
coefficients of homogeneous degree $i$. Since $Y$ is coisotropic, after a linear
change of coordinates $(x_1, \ldots, x_n, y_1, \ldots, y_n)$ we can assume that 
$\alpha_0 = \sum d x_i \wedge d y_i$. 

Note that the ideal $\mathcal{J}$ generated by $(y_1, \ldots, y_q)$ is Poisson with respect to 
the bracket induced by $\alpha$. This means that the coefficients of $d x_r \wedge d x_s$ are in 
$\mathcal{J}$ for each $\alpha_j$ and $r,s \leq q$. We now want to find a formal vector 
field $\mu$ such that the formal diffeomorphism $exp(\mu)$ takes
$\alpha$ to $\alpha_0$ and preserves $\mathcal{J}$. In fact we will construct $exp(\mu)$ 
inductively
as the composition of $exp(\mu_1)$, $exp(\mu_2), \ldots$ where $\mu_i$ is a polynomial vector 
field with coefficients of homogeneous degree $i$ and
$$
exp(\mu_{i-1})\ldots exp(\mu_1) (\alpha) = \alpha_0 + \beta_{\geq i}
$$
with $\beta_{\geq i}$ a 2-form with coefficients of degree $\geq  i$. 
To ensure that $\mathcal{J}$ is preserved, we need to have $\mu_i(\mathcal{J}) 
\subset \mathcal{J}$, which is to say, the coefficient of $\partial/\partial y_r$ in 
$\mu_i$ is an element of $\mathcal{J}$ for $r \leq q$.
Considering $\gamma_i = \alpha_0(\mu_i, \cdot)$ we see that $\gamma_i$ needs to be a polynomial 
differential form with coefficients of degree $i$,
such that the coefficient of $d x_r$ is in $\mathcal{J}$ for $r \leq q$ and $d \gamma_i = \beta_i$
where $\beta_i$ is the degree $=i$ component of $\beta_{\geq i}$. 

Note that $\beta_i$ is at least closed since this property holds for $\alpha$, is preserved 
after the action of $exp(\mu_j)$ and $\beta_i$ is just a homogeneous component of 
the resulting form $\alpha_0 + \beta_{\geq i}$. Since the formal de Rham complex is 
exact in degrees $\geq 0$, $\beta_i = d \gamma_i$ with
 $\gamma_i= \iota_{Eu} \beta_i$, the contraction with the Euler vector field
$Eu = \sum (x_r \partial/\partial x_r + y_r \partial/\partial y_r)$. Note that by induction, after each formal 
diffeomorphism $\mathcal{J}$ remains a Poisson ideal 
hence the coefficient of $d x_r \wedge d x_s$ in $\beta_i$, is an element of $\mathcal{J}$
for $r, s \leq q$. After the Euler field contraction, every coefficient of $d x_r$ in $\gamma_i$ is
also in $\mathcal{J}$, as required. Thus we have a formal diffeomorphism $exp(\mu_i)$ 
that will eliminate $\beta_i$.

Passing to the formal limit $i \to \infty$, we get an 
isomorphism of  $\widehat{\mathcal{O}}_{M, y} \simeq 
k[[x_1, \ldots, x_n, y_1, \ldots y_n]]$  which takes the completion of $I_Y$ to the ideal
$\mathcal{J}$, and is compatible with the symplectic forms.
This proves the ``quasiclassical" part of the statement.

\bigskip
\noindent 
Both $\widehat{\mathcal{O}}_h$ and $\mathcal{D}$ are deformation quantizations
of the same algebra $k[[x_1, \dots, x_n, y_1, \dots, y_n]]$, corresponding to two
formal Poisson bivectors $h (\sum \partial/\partial x_i \wedge \partial/\partial y_i) + h^2\pi_2 + \ldots$
with the same $h$-linear part $\pi_1 = \sum \partial/\partial x_i \wedge \partial/\partial y_i$. 
By the general Maurer-Cartan formalism, Poisson bivectors with 
fixed linear part correspond to Maurer-Cartan solutions of the algebra of polyvector fields 
with the nonzero differential $[\pi_1, \cdot]$. Using $\alpha_0$ to convert polyvector fields to 
differential forms we get the complex in which the bracket with $\pi_1$ becomes the de Rham 
differential. Since the formal de Rham complex is exact in degree two, there is a unique
quantization with the $h$-linear part $\pi_1$. Hence the  above isomorphism modulo $h$ extends to
an isomorphism 
$\sigma_{\mathcal{D}}$. If it can be extended to 
 pair $(\sigma_{\mathcal{D}}, \sigma_{\mathcal{M}})$ compatible with the module
action, then $\widehat{I}_h$, resp, $\mathcal{J}$ is the annihilator of $\widehat{E}_{h, y} / h 
\widehat{E}_{h, y}$, resp. $\mathcal{M}/h \mathcal{M}$, which implies compatibility with 
ideals stated in (1).

\bigskip
\noindent 
For existence of $\sigma_{\mathcal{M}}$ in part (2) assume that compatibility with ideals 
does hold, and  first construct the 
isomorphism modulo $h$  and the lift it inductively modulo higher powers of $h$.
Indeed, using $\sigma_{\mathcal{M}}$ we can view $\widehat{E}_{h, y}/ h \widehat{E}_{h, y}$ and 
	$\mathcal{M}/h
\mathcal{M}$ as projective (hence free) modules of the same rank over
the local ring $\widehat{\mathcal{O}}_{Y, y}$. Hence there is an isomorphism 
$\rho_0:  \mathcal{M}/h
\mathcal{M} \to \widehat{E}_{h, y}/ h \widehat{E}_{h, y} $

To lift it to $\widehat{E}_{h, y}$ and $\mathcal{M}$, take the standard space of generators 
$k^{\oplus e} \subset \mathcal{M}$, and choose
any lift $\rho: k^{\oplus e} \to \widehat{E}_{h, y}$ of $\rho_0|_{k^{\oplus e}}$. 
Consider the subalgebra $\mathcal{D}' \subset \mathcal{D}$ with (topological) generators 
$x_1, \ldots, x_n, y_{q+1}, \ldots, y_n, h$ and the map
$
\sigma': \mathcal{D}' \otimes_k k^{\oplus e} \to \widehat{E}_{h, y}, \quad f \otimes v \mapsto f \cdot \rho(v).
$
By a version of Nakayama's Lemma it is an isomorphism of $k[[h]]$-modules. It induces
an isomorphism with $(\mathcal{D}/\mathcal{D}\langle y_1, \ldots, y_q \rangle)^{\oplus e} = \mathcal{M}$ 
(as $\mathcal{D}$-modules) precisely
when $y_s \cdot \rho(v) = 0$ for any $v \in k^{\oplus e}$ and $s \leq q$. Our goal is to adjust $\rho$ to achieve this condition inductively, ensuring that the vanishing holds modulo $h^{l}$ for 
$l \geq 1$. This obviously works for $l = 1$ as $y_s$ acts by zero on $\widehat{E}_{h, y}/ h \widehat{E}_{h, y}$. To make an inductive step, suppose that we have $\rho_{l+1}: k^{\oplus e} \to \widehat{E}_{h, y}/ h^{l+1} \widehat{E}_{h, y}$ and that by inductive assumption $y_s \cdot \rho_{l+1}(v)$ is divisible by 
$h^l$ for all $v$ and $s \leq q$. Let $U = Im(\rho_{l+1})$ and let 
 $u_1, \ldots, u_e \in U$ be the images of the standard basis vectors. By assumption,
$$
y_s \cdot u_i = h^l \sum_{j} f_{si}^j u_j, \qquad f_{si}^j \in \mathcal{D}'
$$
We are looking for elements 
$$
u'_i = u_i + h^{l-1} \sum_j g_i^j u_j, \qquad g_{i}^j \in \mathcal{D}'
$$ 
which satisfy $y_s \cdot u_i' = 0$ in $\widehat{E}_{h, y}/ h^{l+1} \widehat{E}_{h, y}$. 
For $l \geq 2$ this gives
$$
h f_{si}^j + [y_s, g_{i}^j] = 0\  (mod\ h^2)
$$
Introducing matrices $F_s = (f_{si}^j)$, $G = (g_i^j)$ with entries in $\mathcal{D}'/h \mathcal{D}' 
\simeq k[[x_1, \ldots, x_n, y_{q+1}, \ldots, y_n]]$ and using the fact that $[y_s, \cdot]$ for $s \leq q$
acts as $h (\partial/\partial x_s)$ we get a system of equations $F_s = \partial G/ \partial x_s$ 
on matrices with coefficients in $k[[x_1, \ldots, x_n, y_{q+1}, \ldots, y_n]]$. This has a solution precisely when 
$\partial F_s/ \partial x_t =  \partial F_t /\partial x_s$ for $s, t \leq q$ (by vanishing of formal de Rham 
cohomology in degree 1). But the latter equation is a consequence of $y_t y_s \cdot u_i = y_s y_t \cdot u_i$.
For $l = 1$ the equations read 
$$
h F_s + (y_s G - G y_s) + h G F_s = 0 \  (mod\  h^2)
$$
which is equivalent to $\partial (log(1+G))/\partial x_s = F_s$. As before, we can find a matrix $H$
such that $\partial H/\partial x_s = F_s$ and assume that its entries have zero constant terms. Then
the matrix entries of 
$G = exp(H) -1 $ are formal power series with  zero constant
term. To achieve $y_s \cdot u_i = 0 \  (mod\  h^{l+1})$ we are adjusting $u_i$ by adding vectors that vanish
$(mod\ h^{i-1})$, for all $l \geq 1$. Hence the limit as $l \to \infty$ is well defined, and that gives a system of $\mathcal{D}'$-generators
$u_1, \ldots, u_e$ which also satisfy $y_s \cdot u_i = 0$. This gives an isomorphism  $\sigma_\mathcal{M}$, as
required. 
$\square$

\subsection{Harish-Chandra pairs and torsors  associated to $(\mathcal{O}_h, E_h)$.} 

\bigskip
\noindent 
It is clear from the proofs that isomorphisms of Lemma \ref{isoms} are not unique. Two different choices of the pair
$(\sigma_{\mathcal{D}}, \sigma_{\mathcal{M}})$ are related by automorphisms $\Phi_{\mathcal{D}}: 
\mathcal{D} \to \mathcal{D}, \Phi_{\mathcal{M}}: 
\mathcal{M} \to \mathcal{M}$ compatible with the module action and filtrations. The group $Aut(\mathcal{D}, \mathcal{M})$ formed by all such $(\Phi_{\mathcal{D}}, \Phi_{\mathcal{M}})$ has 
a natural structure of a proalgebraic group, given by reducing the automorphism modulo the image of 
$\mathfrak{m}^k \subset \mathcal{D}$ where $\mathfrak{m}$ is the kernel of the algebra homomorphism 
$\mathcal{D} \to k$ sending $h, x_i, y_j$ to zero. The tangent Lie algebra is formed by pairs $(\phi_{\mathcal{D}}, 
\phi_{\mathcal{M}})$ where $\phi_{\mathcal{D}}: \mathcal{D} \to \mathcal{D}$ is a continuous derivation 
preserving $\mathfrak{m}$ and $\phi_{\mathcal{M}}: \mathcal{M} \to \mathcal{M}$ satisfies 
$$
\phi_{\mathcal{M}} (f m ) = \phi_{\mathcal{D}}(f)  m + f \phi_{\mathcal{M}} (m).
$$
If we drop the condition that $\phi_{\mathcal{D}}$ should preserve $\mathfrak{m}$ we obtain a somewhat larger
Lie algebra $Der(\mathcal{D}, \mathcal{M})$. The following lemma, established by direct
computation, clarifies its structure: 
\begin{lemma} The following short exact sequences hold:
	\begin{enumerate}
		\item 
	The map $(\phi_{\mathcal{D}}, 	\phi_{\mathcal{M}}) \mapsto \phi_{\mathcal{D}}$ induces a short exact sequence
	
\begin{equation}
	\label{DJ}
	0 \to \mathcal{E} = gl_e(\mathcal{D}_p) \to Der(\mathcal{D}, \mathcal{M}) \to Der(\mathcal{D})_{\mathcal{J}} \to 0
\end{equation} 
	where $Der(\mathcal{D})_{\mathcal{J}} $ is the algebra of continuous derivations of $\mathcal{D}$ which 
	sent $\mathcal{J}$ to itself. 
	\item Commutator of  left action of $\phi_{\mathcal{M}}$ with  right action of $\mathcal{E}$ induces the 
	right arrow in the short exact sequence 
	\begin{equation}
	\label{jay}
	0 \to \frac{1}{h} \mathcal{J} \to Der (\mathcal{D}, \mathcal{M}) \to Der(\mathcal{E}) \to 0
	\end{equation}
	\item The reduction of derivations 
	in $Der(\mathcal{E})$, resp.
	in $ Der(\mathcal{D})_{\mathcal{J}}$, modulo $h$, resp. modulo $\mathcal{J}$, 
	induces short exact sequences:
	$$
	0 \to \mathcal{E}/\frac{1}{h} k[[h]] \to Der(\mathcal{E})  \to Ham_q \to 0;
	\qquad 
	0 \to \frac{1}{h} \mathcal{J}/k[[h]] \to Der(\mathcal{D})_{\mathcal{J}} \to Ham_q \to 0
	$$
	where $Ham_q$ is the algebra of Hamiltonian derivations of 
	$k[[x_1, \ldots, x_q, y_1, \ldots, y_q]]$.
 The two compositions $Der(\mathcal{D}, \mathcal{M}) \to Ham_q$ agree and
	this gives rise to a short exact sequence: 
   \begin{equation}
   	\label{product-extension}
   0 \to k[[h]] \to Der(\mathcal{D}, \mathcal{M}) \to Der(\mathcal{D})_{\mathcal{J}} \times_{Ham_q} Der(\mathcal{E}) \to 0 \qquad \square
   \end{equation}

	\end{enumerate}
	
\end{lemma}
Following the pattern of in section 5 in \cite{Ye},  section 6 in \cite{VdB} or 
sections 2, 3 in \cite{DGW} we see that all pairs 
$(\sigma_{\mathcal{D}}, \sigma_{\mathcal{M}})$ are parameterized by a
Harish Chandra torsor $P_{\mathcal{D}, \mathcal{M}}$ over the pair
$(Der(\mathcal{D}, \mathcal{M}), Aut(\mathcal{D}, \mathcal{M}))$. 
Similarly, all isomorphisms $\sigma_{\mathcal{E}}$ are parameterized
by a Harish Chandra torsor $P_{\mathcal{E}}$ over the pair
$(Der(\mathcal{E}), Aut(\mathcal{E}))$. We note here that for 
$P_{\mathcal{D}, \mathcal{M}}$ the connection form $\gamma$ of 
Section 3.2 is an isomorphism (such torsors are called transitive) 
while for $P_{\mathcal{E}}$  the short exact sequence (\ref{jay}) 
implies that this torsor is foliated over $\mathcal{F}$. 

\subsection{The class $\tau_Y$ is the image of $\tau_{Lie}$.}

By Proposition \ref{compare-classes} and Section 4.0.3 in \cite{BNT} the characteristic class 
$$
\tau_Y = \widehat{A}(Q) exp(- \frac{c_1(N)}{2})   e^{-c(\mathcal{O}_h)} ch(E) \in H^\bullet_{DR}(Y)((h)) 
$$
is equal to the image, with respect to the Gelfand-Fuks map of the torsor 
$P_{\mathcal{D}, \mathcal{M}}$, of
the class 
$$
\tau_{Lie} = \widehat{A}_{Lie}(\mathfrak{sp}_{2p})
 exp(- \frac{c_{1, Lie}(\mathfrak{gl}_q)}{2}) e^{-c} 
ch_{Lie}(\mathfrak{gl}_e) 
$$
where the factors other than $e^{-c}$ are defined at the end of Section 3.1 and $c$ is obtained from the extension class
of 
$$
0  \to \frac{1}{h} k[[h]] \to \frac{1}{h} \mathcal{D}\to Der(\mathcal{D})  \to 0
$$
by restricting to the subalgebra $Der(\mathcal{D})_{\mathcal{J}} \subset Der(\mathcal{D}) $ and 
then pulling back under the surjection of (\ref{DJ}). For the factors other than $e^{-c}$ we use the fact
that the components of $\tau_Y$ can be defined by using the Chern-Weil construction on the 
torsor  of symplectic frames in $Q$ and the torsor of usual frames in $N, E$. Since these torsors 
can be induced from $P_{\mathcal{D}, \mathcal{M}}$, we can apply compatibility of Gelfand-Fuks map 
with induced torsors to reinterpret the classes via Lie algebra cohomology of $Der(\mathcal{D}, \mathcal{M})$. 
We record for future reference that 
$$
\tau_{Lie}  \in H^*(Der(\mathcal{D}, \mathcal{M}),  \mathfrak{gl}_q \oplus \mathfrak{gl}_r \oplus \mathfrak{sp}_{2p} 
\oplus \mathfrak{a}'; k((h))).
$$

\subsection{Reduction to class $\tau_{\mathcal{D}_p}$ and end of proof.}

\begin{lemma}
	\label{reduce-class} 
	The cohomology class 
	$
	\tau_{Lie}$ is represented by a
	cocycle which vanishes if one of its arguments is in $\frac{1}{h} \mathcal{J}$. 
  Hence
	$\tau_{Lie}$ is a pullback of a cohomology class of $Der(\mathcal{E})$ via 
	the surjection in (\ref{jay}) and that class
	is further equal to the class of Lemma \ref{vanishing}
	$$
	\tau_{\mathcal{D}_p} \in H^\bullet_{Lie} (Der(\mathcal{E}), \mathfrak{pgl}_r \oplus \mathfrak{sp}_{2p}; k((h)))
	= H^{\bullet}_{Lie} (\g
	, \mathfrak{gl}_r \oplus \mathfrak{sp}_{2p} \oplus \mathfrak{a}'; k((h))).
	$$ 
	\end{lemma}
\textit{Proof.} 

\bigskip
\noindent
\textit{Step 1.} We first recall the definitions. 
Assign the elements in $\mathfrak{gl}_e (k)\subset Der(\mathcal{D}, \mathcal{M})$
 degree 0 and keep assuming that $\deg h = 2$, 
$\deg x_i = \deg y_j = 1$. Then any element of $Der(\mathcal{D}, 
\mathcal{M})$ is a possibly infinite sum of homogeneous elements
of degree $\geq -1$ and the Lie bracket is homogeneous. 

Then 
$\mathfrak{gl}_e$ and $\mathfrak{sp}_{2p}$ are spanned by 
elements $\frac{1}{h} x_i y_j$, $1 \leq i, j \leq q$ and 
$\frac{1}{h} x_s x_t, \frac{1}{h} y_s y_t, \frac{1}{h} (x_s y_t + y_t x_s)$, 
$ q+1 \leq s, t \leq p+q$, respectively, and the degree zero part 
of $Der(\mathcal{D}, \mathcal{M})$ splits as
$$
\mathfrak{gl}_e \oplus \mathfrak{gl}_q \oplus \mathfrak{sp}_{2p} 
\oplus W
$$
where $W$ is spanned by $\frac{1}{h} y_j y_t, \frac{1}{h} x_s y_j $
with $1 \leq j \leq q$, $1 \leq t \leq (p+q)$, $ (q+1) \leq s \leq (p+q)$
(we note here that in the specified ranges the variables commute). This gives
a projection 
$$
Der(\mathcal{D}, \mathcal{M}) \to \mathfrak{gl}_e \oplus \mathfrak{gl}_q \oplus \mathfrak{sp}_{2p}
$$ 
sending the elements of nonzero degree to zero, and vanishing on $W$. We can 
combine it with a natural projection to any of the three factors on the right hand
side, to be used for calculation of classes 
$ ch_{Lie}(\mathfrak{gl}_e)$,  $\widehat{A}_{Lie}(\mathfrak{sp}_{2p})$, 
$exp(- \frac{c_{1, Lie}(\mathfrak{gl}_q)}{2})$ in the definition of $\tau_{Lie}$.

The curvature defined in (\ref{curvature})  
 is not zero only when its arguments have degrees $-1, 0$ or $1$. The same applies to the 
degree zero component $c_0$ of $c$. We recall here that 
 $c_0$ is computed with respect to the projection onto $k$ which vanishes
on elements of non-zero degrees, trace zero matrices in $\mathfrak{gl}_e$ and on the subspaces $\mathfrak{sp}_{2p}$, $W$, 
and on the elements of the type $\frac{1}{h} (x_i y_j + y_j x_i)$.

\bigskip
\noindent 
\textit{Step 2.} Let us show that the $\widehat{A}_{Lie}$ class is pulled back 
from the quotient by $\frac{1}{h} \mathcal{J}$.  In fact, consider 
the curvature $C(u \wedge v) \in \h$ for the projection onto $\h = 
\mathfrak{sp}_{2p}$ 
and $u \in \frac{1}{h} \mathcal{J}$ of degree $-1, 0$ or $1$. 
It follows from (\ref{jay}) that $\frac{1}{h} \mathcal{J}$ is a Lie ideal which has 
zero projection onto $\mathfrak{sp}_{2p}$, so  all positive components of $\widehat{A}_{Lie}(\mathfrak{sp}_{2p})$
vanish if one of the arguments is in $\frac{1}{h} \mathcal{J}$.

\bigskip
\noindent
\textit{Step 3.} Let us prove that the cochain representing the class 
$$
ch_{Lie}(\mathfrak{gl}_e) exp(- \frac{1}{2}c_{1, Lie}(\mathfrak{gl}_q) - c_0) 
$$
is zero if one of its arguments is in $\frac{1}{h} \mathcal{J}$. Since 
$c_0$ vanishes on elements $\frac{1}{h} (x_i y_j + y_j x_i)$
we have
$$
c_0 (\frac{1}{h}\sum b_{ij} x_i y_j) = 
c_0(\frac{1}{2h} \sum b_{ij} (x_i y_j - y_j x_i) + \frac{1}{2h} \sum b_{ij} (x_i y_j + y_j x_i)) = - \frac{1}{2} \sum b_{ii} 
$$
Since $\sum b_{ii}$ is the invariant polynomial corresponding to 
$c_{1, Lie} (\mathfrak{gl}_q)$ we conclude that 
$c_0 + c_{1, Lie} (\mathfrak{gl}_q)$ corresponds to 
the linear function on $\h = \mathfrak{gl}_e \oplus \mathfrak{gl}_q$ which sends
$(X_1, X_2)$ to $\frac{1}{e} tr(X_1)$. 
Moreover, since $ch_{Lie}(\mathfrak{gl}_e)$ comes from $tr(exp(x))$ and 
$tr(exp(x - \alpha \cdot I)) = tr(exp(x)) exp (-\alpha)$, we can rewrite the above class as 
the image of the invariant series 
$$
S(X_1 \oplus X_2) = \sum_{n \geq 0} \frac{1}{n!} tr \big(X_1 - \frac{1}{e} tr(X_1)\big)^n
$$
under the Chern-Weil map. We denote by $\overline{X}
= X_1 - \frac{1}{e} tr(X_1)$ the trace zero part of $X_1$ and by $S_l(X) 
= \frac{1}{l!} tr (\overline{X}_1)^l$ the
degree $l$ component of the invariant power series. 
Recall that the Chern-Weil class corresponding to $S_l(X)$ is obtained by polarization of
$S_l(X)$:
$$
\rho(S_l)(v_1 \wedge \ldots \wedge v_{2l} )= 
\frac{1}{l!} \sum_{\sigma} (-1)^\sigma  
tr( \overline{C}(v_{\sigma(1)} \wedge v_{\sigma(2)})\overline{C}(( v_{\sigma(3)}\wedge v_{\sigma(4)})\ldots\overline{C}(v_{\sigma(2l-1)}\wedge v_{\sigma(2l)}))
$$
where the sum is over all permutations $\sigma \in S_{2n}$ that 
satisfy $\sigma(2i-1) < \sigma(2i)$. So it suffices to show that 
$\overline{C} = 0$ in $\mathfrak{gl}_e$ if 
$C = C(u \wedge v)$ with $u \in \frac{1}{h} \mathcal{J}$. This holds since 
$\frac{1}{h} \mathcal{J}$ is a Lie ideal and its projection 
onto $\mathfrak{gl}_q$ lands into the subspace of scalar matrices which have
trivial $\overline{X}$ part. 

\bigskip
\noindent 
\textit{Step 4.} It remains to show that the class $exp(-(c - c_0))$ is
in the image of the pullback under the projection 
$Der(\mathcal{D}, \mathcal{M}) \to Der(\mathcal{E})$. 
Recall that $c$ was defined in Section 4.3.
Let $\theta$ be the pullback of a similar class 
under $Der(\mathcal{D}, \mathcal{M}) \to Der(\mathcal{E})$. Then by the 
short exact sequence (\ref{product-extension}), the sum $c - \theta$ is zero 
(the minus sign is in front of $\theta$ is due to the fact that $\mathcal{M}$ was
considered as a right $\mathcal{E}$-module or, equivalently, a left 
$\mathcal{E}^{op}$-module). We are using the fact that the sum of two extensions 
descends to the fiber product over $Ham_p$, that the the pullback of the extension class 
to the extension algebra $Der(\mathcal{D}, \mathcal{M})$ must be zero, and that 
existence of $P_{\mathcal{D}, \mathcal{M}}$ allows to use the Gelfand-Fuks map 
associated to this torsor. Hence, on the level of cohomology
$c = \theta$ and the same holds for each of the coeffiients in the expantion in 
powers of $h$, e.g. $c_0 = \theta_0$.

The fact that the class pulled back from $Der(\mathcal{E})$ is exactly $\tau_{\mathcal{D}_p}$ 
follows from the definition of $\tau_{Lie}$ and the vanishing proved. 
$\square$ 

\bigskip
\noindent
\textbf{End of proof of Theorem \ref{main}.} By Section 2.2 (application of 
Riemann Roch Theorm), the class $\tau(E_h)$ is the 
image of $\tau_Y(E)$ in the cohomology of $Y$.

By Proposition \ref{compare-classes} the class $\tau_Y(E)$ in the 
de Rham cohomology of $Y$ is the image of a class $\tau_{Lie}$ under
the Gelfand-Fuks map associated to the Harish-Chandra torsor $P_{\mathcal{D}, \mathcal{M}}$ of 
all isomorphisms $(\sigma_{\mathcal{D}}, \sigma_{\mathcal{M}})$.

By Lemma \ref{reduce-class}, we can replace the pair $(\tau_{Lie}, P_{\mathcal{D}, \mathcal{M}})$ by the 
pair $(\tau_{\mathcal{D}_p}, P_{\mathcal{E}})$ where the torsor $P_{\mathcal{E}}$ parameterizes isomorphisms $\sigma_{\mathcal{E}}:
	End_{\widehat{\mathcal{O}}_h} 
	\widehat{E}_h \to \mathcal{E}$ (we note here that a choice of $(\sigma_{\mathcal{D}}, \sigma_{\mathcal{M}})$ also induces
	a choice of $\sigma_{\mathcal{E}}$). 
	
	But $P_{\mathcal{E}}$  is a Harish-Chandra torsor foliated over
	$\mathcal{F} \subset T_Y$, hence the characteristic clas $\tau_Y$ admits a lift to 
	$\bigoplus_{r \geq 0}^p H^{2r}(Y, F^r \Omega^\bullet_Y)  ((h))$
	by  Proposition \ref{compare-classes}. 
	
	Finally, the components of $\tau_{\mathcal{D}_p}$ vanish in cohomology groups 
	of degrees $> 2p$ by 
	Lemma \ref{vanishing}. This finishes the proof of Theorem \ref{main} $\square$

\bigskip
\noindent
Email: vbaranov@math.uci.edu
\end{document}